\documentclass{article}
\begin{document}

\begin{center}
{\textbf{ON GENERALIZED ENTROPY MEASURES AND PATHWAYS}}

\medskip

A.M. MATHAI

\smallskip

Department of Mathematics and Statistics, McGill University, Montreal, Canada H3A 2K6, and \\
Centre for Mathematical Sciences, Pala Campus, Arunapuram P.O., Pala-686 574, Kerala, India

\smallskip

and

\smallskip
H.J. HAUBOLD

\smallskip

Office for Outer Space Affairs, United Nations, Vienna International Centre, P.O. Box 500, A-1400 Vienna, Austria and\\
Centre for Mathematical Sciences, Pala Campus, Arunapuram P.O., Pala-686 574, Kerala, India

\smallskip

\end{center}

\noindent
\textbf{Abstract.} Product probability property, known in the literature as 
statistical independence, is examined first. Then  
generalized entropies are introduced, all of which give 
generalizations to Shannon entropy. It is shown that the nature 
of the recursivity postulate automatically determines the 
logarithmic functional form for Shannon entropy. Due to the 
logarithmic nature, Shannon entropy naturally gives rise to 
additivity, when applied to situations having product 
probability property. It is argued that the natural process is 
non-additivity, important, for example, in statistical mechanics (Tsallis 2004, Cohen 2005), 
even in product probability property situations and additivity
can hold due to the involvement of a recursivity postulate 
leading to a logarithmic function. Generalized entropies 
are introduced and some of
their properties are examined. Situations are examined 
where a generalized entropy of order $\alpha$ leads to pathway models, exponential and power 
law behavior and related differential 
equations. Connection of this entropy to Kerridge's measure 
of ``inaccuracy" is also explored.\\

\medskip


\noindent
\textbf{1.\hskip.4cm Introduction}

\medskip

\noindent

Mathai and Rathie (1975) consider various 
generalizations of Shannon entropy (Shannon, 1948), called 
entropies of order $\alpha$, and give various properties, 
including additivity property, and characterization theorems. 
Recently, Mathai and Haubold (2006, 2006a) 
explored a generalized entropy of 
order $\alpha$, which is connected to a measure of uncertainty 
in a probability scheme, Kerridge's (Kerridge, 1961) concept of inaccuracy
in a scheme, and pathway models that are considered in this paper.

\medskip

As defined in Mathai and Haubold (2006, 2006a) the entropy $M_{k,\alpha}(P)$ is a non-additive entropy and
his measure $M_{k,\alpha}^{*}(P)$ is an additive entropy. It is
also shown that maximization of the continuous analogue of 
$M_{k,\alpha}(P)$, denoted by $M{\alpha}(f)$, gives rise to 
various functional forms for $f$, depending upon the types of 
constraints on $f$.

\medskip
 
Occasionally, emphasis is placed on the fact that Shannon entropy satisfies
the additivity property, leading to extensivity. It will be shown that
when the product probability property (PPP) holds then a logarithmic
function can give a sum and a logarithmic function enters into Shannon
entropy due to the assumption introduced through a certain type of
recursivity postulate. The concept of statistical independence will be
examined in Section 1 to illustrate that simply because of PPP one
need not expect additivity to hold or that one should not expect this
PPP should lead to extensivity. The types of non-extensivity,
associated with a number of generalized entropies, are pointed out even
when PPP holds. The nature of non-extensivity that can be expected
from a multivariate distribution, when PPP holds or when there is
statistical independence of the random variables, is illustrated by
taking a trivariate case.  

\medskip
 
Maximum entropy principle is examined in Section 2. It is 
shown that optimization of measures of entropies, in 
the continuous populations, under selected constraints,
leads to various types of models. It is shown that the generalized 
entropy of order $\alpha$ is a convenient one to obtain various 
probability models.
 
\medskip
 
Section 3 examines the types of differential equations 
satisfied by the various special cases of the pathway model.
\medskip

\noindent
\textbf{1.1.\hskip.4cm Product probability property (PPP) or 
statistical independence of events}
\medskip

Let $P(A)$ denote the probability of the event $A$. If the 
definition $P(A\cap B)=P(A)P(B)$ is taken as the definition of 
independence of the events $A$ and $B$ then any event $A\in S$,
and $S$ the sure event are independent. But $A$ is contained 
in $S$ and then the definition of independence becomes 
inconsistent with the common man's vision of independence. Even
if the trivial cases of the sure event $S$ and the impossible 
event $\phi$ are deleted, still this definition becomes a 
resultant of some properties of positive numbers. Consider a 
sample space of $n$ distinct elementary events. If symmetry in 
the outcomes is assumed then we will assign equal probabilities 
$\frac{1}{n}$ each to the elementary events. Let $C=A\cap B$. 
If $A$ and $B$ are independent then $P(C)=P(A)P(B)$. Let
$$
P(A)=\frac{x}{n},~P(B)=\frac{y}{n},~P(C)=\frac{z}{n}.
$$
Then
\begin{equation}
\left(\frac{x}{n}\right)\left(\frac{y}{n}\right)
=\left(\frac{z}{n}\right)\Rightarrow nz=xy,~x,y,z=1,2,...,n-1,
z<x,y
\end{equation}
deleting $S$ and $\phi$. There is no solution for $x,y,z$ for a
large number of $n$, for example, $n=3,5,7$. This means that 
there are no independent events in such cases and it sounds 
strange from a common man's point of view.
\medskip

The term ``independence'' of events is a misnomer. This property
should have been called product probability property or PPP of 
events. There is no reason to expect the 
information or entropy in a joint distribution to be the sum of
the information contents of the marginal distributions when the
PPP holds for the distributions, that is when the joint density or 
probability function is a product of the marginal densities or 
probability functions. We may expect a term due to the product 
probability to enter into the expression for the entropy in the
joint distribution in such cases. But if the information or 
entropy is defined in terms of a logarithm, then naturally, 
logarithm of a product being the sum of logarithms, we can expect
a sum coming in such situations. This is not due to independence 
or due to the PPP of the densities but due to the fact that a 
functional involving logarithm is taken thereby a product has 
become a sum. Hence not too much importance should be put on 
whether or not the entropy on the joint distribution becomes 
sum of the entropies on marginal distributions or additivity 
property when PPP holds.
\medskip

\noindent
\textbf{1.2.\hskip.4cm How is logarithm coming in Shannon's 
entropy?}
\medskip

Several characterization theorems for Shannon entropy and its 
various generalizations are given in Mathai and Rathie (1975. 
Modified and refined versions of Shannon's own postulates
are given as postulates for the first theorem characterizing 
Shannon entropy in Mathai and Rathie (1975). Apart from 
continuity, symmetry, zero-indifference and normalization 
postulates the main postulate in the theorem is a recursivity 
postulate, which in essence says that when the PPP holds then 
the entropy will be a weighted sum of the entropies, thus in 
effect, assuming a logarithmic functional form. The crucial 
postulate is stated here. Consider a multinomial population 
$P=(p_1,...,p_m)$, $p_i>0,~i=1,...,m$, $p_1+...+p_m=1$, 
that is, $p_i=P(A_i),~i=1,...,m$, $A_1\cup...\cup A_m=S$, 
$A_i\cap A_j= \phi,~i\ne j$. If any $p_i$ can take a zero value
also then zero-indifferent postulate, namely that the entropy 
remains the same when an impossible event is incorporated into 
the scheme, is to be added. Let $H_n(p_1,...,p_n)$ denote the 
entropy to be defined. Then the crucial recursivity postulate 
says that
\begin{eqnarray}
&&H_n(p_1,...,p_{m-1},p_mq_1,..,p_mq_{n-m+1})\nonumber\\
&=&H_m(p_1,...,p_m)
+p_mH_{n-m+1}(q_1,...,q_{n-m+1})
\end{eqnarray}
$\sum_{i=1}^mp_i=1,~\sum_{i=1}^{n-m+1}q_i=1$. This says that if
the $m$-th event $A_m$ is partitioned into independent events 
$P(A_m\cap B_j)=P(A_m)P(B_j)=p_mq_j,~j=1,...,n-m+1$ so that 
$p_m=p_mq_1+...+p_mq_{n-m+1}$ then the entropy $H_n(\cdot)$ 
becomes a weighted sum. Naturally, the result will be a 
logarithmic function for the measure of entropy.

\medskip

There are several modifications to this crucial recursivity 
postulate. One suggested by Tverberg is that $n-m+1=2$ and 
$q_1=q, q_2=1-q, ~0<q<1$ and $H_2(q,1-q)$ is assumed to be 
Lebesgue integrable in $0\le q\le 1$. Again a characterization 
of Shannon entropy is obtained. In all the characterization 
theorems for Shannon entropy this recursivity property enters 
in one form or the other as a postulate, which in effect implies
a logarithmic form for the entropy measure. Shannon entropy 
$S_k$ has the following form:
\begin{equation}
S_k=-A\sum_{i=1}^kp_i\ln p_i, ~p_i>0,~i=1,...,k, p_1+...+p_k=1,
\end{equation}
where $A$ is a constant. If any $p_i$ is assumed to be zero 
then $0\ln 0$ is to be interpreted as zero. Since the constant 
$A$ is present, logarithm can be taken to any base. Usually the
logarithm is taken to the base $2$ for ready application to 
binary systems. We will take logarithm to the base ${\rm e}$.
\medskip

\noindent
\textbf{1.3.\hskip.4cm Generalization of Shannon entropy}
\medskip

Consider again a multinomial population $P=(p_1,...,p_k),~p_i>0,
i=1,...,k, p_1+...+p_k=1$. The following are some of the 
generalizations of Shannon entropy $S_k$.
\begin{eqnarray}
R_{k,\alpha}(P)&=&\frac{\ln(\sum_{i=1}^kp_i^{\alpha})}{1-\alpha},
~\alpha\ne 1, \alpha>0, \\
&&\mbox{  (R\'enyi entropy of order $\alpha$ of 1961)}\nonumber\\
H_{k,\alpha}(P)&=&\frac{\sum_{i=1}^kp_i^{\alpha}-1}{2^{1-\alpha}-1},
~\alpha\ne 1,~\alpha>0\\
& &\nonumber\;\;\mbox{  (Havrda-Charv\'at  entropy of order $\alpha$ of 1967)}\\T_{k,\alpha}(P)&=&\frac{\sum_{i=1}^kp_i^{\alpha}-1}{1-\alpha},
~\alpha\ne 1, ~\alpha>0\\
&&\mbox{  (Tsallis entropy of 1988)}\nonumber\\
M_{k,\alpha}(P)&=&\frac{\sum_{i=1}^kp_i^{2-\alpha}-1}{\alpha-1},
~\alpha\ne 1,~-\infty<\alpha<2\\
&& \mbox{  (entropic form of order $\alpha$)}\nonumber\\
M_{k,\alpha}^{*}(P)&=&\frac{\ln(\sum_{i=1}^kp_i^{2-\alpha})}{\alpha-1},
~\alpha\ne 1, ~-\infty<\alpha<2, \\
&&\;\;\mbox{  (additive entropic form of order $\alpha$).}\nonumber
\end{eqnarray}
When $\alpha\rightarrow 1$ all the entropies of order $\alpha$ 
described above in (4) to (7) go to Shannon entropy $S_k$.
\begin{equation}
\lim_{\alpha\rightarrow 1}R_{k,\alpha}(P)
=\lim_{\alpha\rightarrow 1}H_{k,\alpha}(P)
=\lim_{\alpha\rightarrow 1}T_{k,\alpha}(P)
=\lim_{\alpha\rightarrow 1}M_{k,\alpha}(P)
=\lim_{\alpha\rightarrow 1}M_{k,\alpha}^{*}(P)=S_k.
\end{equation}
Hence all the above measures are called generalized entropies 
of order $\alpha$.
\medskip

Let us examine to see what happens to the above entropies in 
the case of a joint distribution. Let $p_{ij}>0,~i=1,...,m, 
j=1,...,n$ such that $\sum_{i=1}^m\sum_{j=1}^np_{ij}=1$. This 
is a bivariate situation of a discrete distribution. Then the 
entropy in the joint distribution, for example,
\begin{equation}
M_{m,n,\alpha}(P,Q)=\frac{\sum_{i=1}^m\sum_{j=1}^np_{ij}^{2-\alpha}-1}{\alpha-1}.
\end{equation}
If the PPP holds and if $p_{ij}=p_iq_j$, $p_1+...+p_m=1$, 
$q_1+...+q_n=1$, $p_i>0,~i=1,...,m$, $q_j>0, ~j=1,...,n$ and 
if $P=(p_1,...,p_m), ~Q=(q_1,...,q_n)$ then
\begin{eqnarray*}
(\alpha-1)M_{m,\alpha}&(P)& M_{n,\alpha}(Q)=\frac{1}{\alpha-1}
\left(\sum_{i=1}^mp_i^{2-\alpha}-1\right)
\left(\sum_{j=1}^nq_j^{2-\alpha}-1\right)\\
&=&\frac{1}{\alpha-1}\left[\sum_{i=1}^m\sum_{j=1}^np_i^{2-\alpha}q_j^{2-\alpha}
-\sum_{i=1}^mp_i^{2-\alpha}-\sum_{j=1}^nq_j^{2-\alpha}+1\right]\\
&=&M_{m,n,\alpha}(P,Q)-M_{m,\alpha}(P)-M_{n,\alpha}(Q).
\end{eqnarray*}
Therefore
\begin{equation}
M_{m,n,\alpha}(P,Q)=M_{m,\alpha}(P)+M_{n,\alpha}(Q)
+(\alpha-1)M_{m,\alpha}(P)M_{n,\alpha}(Q).
\end{equation}
If any one of the above mentioned generalized entropies in 
(4) to (8) is written as $F_{m,n,\alpha}(P,Q)$ then we 
have the relation
\begin{equation}
F_{m,n,\alpha}(P,Q)=F_{m,\alpha}(P)+F_{n,\alpha}(Q)
+a(\alpha)F_{m,\alpha}(P)F_{n,\alpha}(Q).
\end{equation}
where
\begin{eqnarray}
a(\alpha)&=&0\; \mbox{(R\'enyi entropy $R_{k,\alpha}(P)$)}\nonumber\\
&=&2^{1-\alpha}-1\ \ \mbox{(Havrda-Charv\'at entropy $H_{k,\alpha}(P)$)}\nonumber\\
&=&1-\alpha\ \ \mbox{(Tsallis entropy $T_{k,\alpha}(P)$)}\nonumber\\
&=&\alpha-1\ \ \mbox{(entropic form of order $\alpha$, i.e., $M_{k,\alpha}(P)$)}\nonumber\\
&=&0\ \ \mbox{(additive entropic form of order $\alpha$, i.e., $M_{k,\alpha}^{*}(P)$)}.
\end{eqnarray}
When $a(\alpha)=0$ the entropy is called additive and when 
$a(\alpha)\ne 0$ the entropy is called non-additive. As can be 
expected, when a logarithmic function is involved, as in the 
cases of $S_k(P),R_{k,\alpha}(P), M_{k,\alpha}^{*}(P)$, the 
entropy is additive and $a(\alpha)=0$.
\medskip

\noindent
\textbf{1.4.\hskip.4cm Extensions to higher dimensional joint 
distributions}
\medskip

Consider a trivariate population or a trivariate discrete 
distribution $p_{ijk}>0$,$~i=1,...,m$,$~j=1,...,n$,$~k=1,...,r$ 
such that $\sum_{i=1}^m\sum_{j=1}^n\sum_{k=1}^rp_{ijk}=1$. If 
the PPP holds mutually, that is, pair-wise as well as jointly, 
which then will imply that
\begin{eqnarray*}
p_{ijk}&=&p_iq_js_k,~\sum_{i=1}^mp_i=1,~\sum_{j=1}^nq_j=1,
~\sum_{k=1}^rs_k=1,\\
P&=&(p_1,...,p_m),~Q=(q_1,...,q_n),~ S=(s_1,...,s_r).
\end{eqnarray*}
Then proceeding as before, we have for any of the measures 
described above in (4) to (8), calling it $F(\cdot)$,
\begin{eqnarray}
F_{m,n,r,\alpha}(P,Q,S)&=&F_{m,\alpha}(P)+F_{n,\alpha}(Q)
+F_{r,\alpha}(S)+a(\alpha)[F_{m,\alpha}(P)F_{n,\alpha}(Q)\nonumber\\
&& +F_{m,\alpha}(P)F_{r,\alpha}(S)
+F_{n,\alpha}(Q)F_{r,\alpha}(S)]\nonumber\\
&&+[a(\alpha)]^2F_{m,\alpha}(P)F_{n,\alpha}(Q)F_{r,\alpha}(S)
\end{eqnarray}
where $a(\alpha)$ is the same as in (13). The same procedure 
can be extended to any multivariable situation. If $a(\alpha)=0$
we may call the entropy additive and if $a(\alpha)\ne 0$ then 
the entropy is non-additive.
\medskip

\noindent
\textbf{1.5.\hskip.4cm Crucial recursivity postulate}
\medskip

Consider the multinomial population $P=(p_1,...,p_k),~p_i>0,
~i=1,...,k,~p_1+...+p_k=1$. Let the entropy measure to be 
determined through appropriate postulates be denoted by 
$H_k(P)=H_k(p_1,...,p_k)$. For $k=2$ let
\begin{equation}
f(x)=H_2(x,1-x),~0\le x\le 1\mbox{  or  } x\in [0,1].
\end{equation}
If another parameter $\alpha$ is to be involved in $H_2(x,1-x)$
then we will denote $f(x)$ by $f_{\alpha}(x)$.  From (5) to 
(7) it can be seen that the generalized entropies  of order 
$\alpha$ of Havrda-Charv\'at (1967), Tsallis (1988, 2004) and Shannon (1948) 
entropy satisfy the functional equation
\begin{equation}
f_{\alpha}(x)+b_{\alpha}(x)f_{\alpha}\left(\frac{y}{1-x}\right)
=f_{\alpha}(y)+b_{\alpha}(x)f\left(\frac{x}{1-y}\right)
\end{equation} 
for $x,y\in [0,)$ with $x+y\in [0,1]$, with the boundary 
condition
\begin{equation}
f_{\alpha}(0)=f_{\alpha}(1)
\end{equation}
where
\begin{eqnarray}
b_{\alpha}(x)&=&1-x\mbox{  (Shannon entropy $S_k(P)$)}\nonumber\\
&=&(1-x)^{\alpha}\mbox{  (Harvda-Charv\'at entropy $H_{k,\alpha}(P)$)}\nonumber\\
&=&(1-x)^{\alpha}\mbox{  (Tsallis entropy $T_{k,\alpha}(P)$)}\nonumber\\
&=&(1-x)^{2-\alpha}\mbox{  (entropic form of order $\alpha$, i.e., $M_{k,\alpha}(P)$)}.
\end{eqnarray}
Observe that the normalizing constant at $x=\frac{1}{2}$ is 
equal to $1$ for $H_{k,\alpha}(P)$ and it is different for 
other entropies. Thus equations (6),(7),(8), with the 
appropriate normalizing constants $f_{\alpha}(\frac{1}{2})$, 
can give characterization theorems for the various entropy 
measures. The form of $b_{\alpha}(x)$ is coming from the crucial
 recursivity postulate, assumed as a desirable property for the 
measures.
\medskip

\noindent
\textbf{1.6.\hskip.4cm Continuous analogues}
\medskip

In the continuous case let $f(x)$ be the density function of a 
real random variable $x$. Then the various entropy measures, 
corresponding to the ones in (4) to (8) are the following:
\begin{eqnarray}
R_{\alpha}(f)&=&\frac{1}{1-\alpha}\ln\left[\int_{-\infty}^{\infty}[f(x)]^{\alpha}{\rm d}x\right],
~\alpha\ne 1,~ \alpha>0\\
&&\mbox{  (R\'enyi entropy of order $\alpha$)}\nonumber\\
H_{\alpha}(f)&=&\frac{1}{2^{1-\alpha}-1}\left[\int_{-\infty}^{\infty}[f(x)]^{\alpha}{\rm d}x-1\right],
~\alpha\ne 1, ~\alpha>0\\
&&\;\;\nonumber\mbox{  (Havrda-Charv\'at entropy of order $\alpha$)}\\
T_{\alpha}(f)&=&\frac{1}{1-\alpha}\left[\int_{-\infty}^{\infty}[f(x)]^{\alpha}{\rm d}x-1\right],
~\alpha\ne 1, ~\alpha>0,\\
&&\mbox{  (Tsallis entropy of order $\alpha$)}\nonumber\\
M_{\alpha}(f)&=&\frac{1}{\alpha-1}\left[\int_{-\infty}^{\infty}[f(x)]^{2-\alpha}{\rm d}x-1\right],
~\alpha\ne 1,~\alpha<2\\
&&\mbox{  (entropic form of order $\alpha$)}\nonumber\\
M_{\alpha}^{*}(f)&=&\frac{1}{\alpha-1}\ln\left[\int_{-\infty}^{\infty}[f(x)]^{2-\alpha}{\rm d}x\right],
~\alpha\ne 1, ~\alpha<2\\
&&\nonumber\;\;\mbox{  (additive entropic form of order $\alpha$)}.
\end{eqnarray}
As expected, Shannon entropy in this case is given by
\begin{equation}
S(f)=-A\int_{-\infty}^{\infty}f(x)\ln f(x){\rm d}x
\end{equation}
where $A$ is a constant.\\

Note that when PPP (product probability property) or statistical
independence holds then in the continuous case also we have the 
property in (12) and (14) and then non-additivity holds for 
the measures analogous to the ones in (3),(5),(6),(7)
with $a(\alpha)$ remaining the same. Since the steps are 
parallel a separate derivation is not given here.
\medskip

\noindent
\textbf{2.\hskip.4cm Maximum Entropy Principle}
\medskip

If we have a multinomial population $P=(p_1,...,p_k),~p_i>0,
i=1,...,k,p_1+...+p_k=1$ or the scheme $P(A_i)=p_i$, 
$A_1\cup...\cup A_k=S$, $P(S)=1$, $A_i\cap A_j= \phi,~i\ne j$ 
then we know that the maximum uncertainty in the scheme or the 
minimum information from the scheme is obtained when we cannot 
give any preference to the occurrence of any particular event 
or when the events are equally likely or when $p_1=p_2=...=p_k
=\frac{1}{k}$. In this case, Shannon entropy becomes,
\begin{equation}
S_k(P)=S_k(\frac{1}{k},...,\frac{1}{k})=-A\sum_{i=1}^k\frac{1}{k}\ln \frac{1}{k}=A\ln k
\end{equation}
and this is the maximum uncertainty or maximum Shannon entropy 
in this scheme. If the arbitrary functional $f$ is to be fixed 
by maximizing the entropy then in (19) to (21) we have to 
optimize $\int_{-\infty}^{\infty}[f(x)]^{\alpha}{\rm d}x$ for 
fixed $\alpha$, over all functional $f$, subject to the condition
$\int_{-\infty}^{\infty}f(x){\rm d}x=1$ and $f(x)\ge 0$ for all
$x$. For applying calculus of variation procedure we consider 
the functional

$$
U=[f(x)]^{\alpha}-\lambda[f(x)]
$$
where $\lambda$ is a Lagrangian multiplier. Then the Euler 
equation is the following:
\begin{equation}
\frac{\partial U}{\partial f}=0\Rightarrow
\alpha f^{\alpha-1}-\lambda=0\Rightarrow 
f=\left(\frac{\lambda}{\alpha}\right)^{\frac{1}{\alpha-1}}
=\mbox{ constant}.
\end{equation}
Hence $f$ is the uniform density in this case, analogous to the 
equally likely situation in the multinomial case. If the first 
moment $E(x)=\int_{-\infty}^{\infty}xf(x){\rm d}x$ is assumed to
be a given quantity for all functional $f$ then $U$ will become
the following for (19) to (21).

$$
U=[f(x)]^{\alpha}-\lambda_1[f(x)]-\lambda_2 xf(x)
$$
and the Euler equation leads to the power law. That is,
\begin{equation}
\frac{\partial U}{\partial f}=0\Rightarrow 
\alpha f^{\alpha-1}-\lambda_1-\lambda_2 x=0\Rightarrow 
f=c_1\left[1+\frac{\lambda_2}{\lambda_1}x\right]^{\frac{1}{\alpha-1}}.
\end{equation}
By selecting $c_1,\lambda_1,\lambda_2$ appropriately we can 
create a density out of  (27). For $\alpha>1$ and 
$\frac{\lambda_2}{\lambda_1}>0$ the right side in (27) 
increases exponentially. If $\alpha=q>1$ and 
$\frac{\lambda_2}{\lambda_1}=q-1$ then we have Tsallis' 
$q$-exponential function from the right side of (27). 
If $\alpha>1$ and $\frac{\lambda_2}{\lambda_1}=-(\alpha-1)$ 
then (27) can produce a density in the category of a type-1 
beta. From (27) it is seen that the form of the entropies of  
Havrda-Charv\'at $H_{k,\alpha}(P)$ and Tsallis $T_{k,\alpha}(P)$ 
need special attention to produce densities (Ferri et al. 2005). However, Tsallis  
has considered a different constraint on $E(x)$. 
If the density $f(x)$ is replaced by its escort density, namely, 
$\mu[f(x)]^{\alpha}$ where $\mu^{-1}=\int_x[f(x)]^{\alpha}{\rm d}x$ 
and if the expected value of $x$ in this escort density is 
assumed to be fixed for all functional $f$ then the $U$ of 
(26) becomes
\begin{eqnarray*}
U&=f^{\alpha}-\lambda_1 f+\mu \lambda_2 xf^{\alpha}\\
\mbox{and}\\
\frac{\partial U}{\partial f}&=0\Rightarrow
\alpha f^{\alpha-1}[1+\mu\lambda_2 x]=\lambda_1\Rightarrow 
f=\frac{\lambda^{*}}{(1+\lambda_3 x)^{\frac{1}{\alpha-1}}}
\Rightarrow\\
f&=\lambda_1^{*}[1+\lambda_3 x]^{-\frac{1}{\alpha-1}}
\end{eqnarray*}
where $\lambda_3$ is a constant and $\lambda_1^{*}$ is the 
normalizing constant. If $\lambda_3$ is taken as 
$\lambda_3=\alpha-1$ then
\begin{equation}
f=\lambda_1^{*}[1+(\alpha-1)x]^{-\frac{1}{\alpha-1}}.
\end{equation}
Then (28) for $\alpha>1$ is Tsallis statistics (Tsallis 2004, Cohen 2005). 
Then for $\alpha<1$ also by writing $\alpha-1=-(1-\alpha)$ one 
gets the case of Tsallis statistics for $\alpha<1$ (Ferri et al. 2005). These 
modifications and the consideration of escort distribution are 
not necessary if we take the generalized entropy of order $\alpha$. Thus if we 
consider $M_{\alpha}(f)$ and if we assume that the first moment 
in $f(x)$ itself is fixed for all functional $f$ then the Euler 
equation gives

$$
(2-\alpha)f^{1-\alpha}-\lambda_1+\lambda_2 x=0\Rightarrow 
f=\bar{\lambda}\left[1-\frac{\lambda_2}{\lambda_1}x\right]^{\frac{1}{1-\alpha}}
$$
and for $\frac{\lambda_2}{\lambda_1}=1-\alpha$ we have Tsallis 
statistics (Tsallis 2004, Cohen 2005)
\begin{equation}
f=\bar{\lambda}[1-(1-\alpha)x]^{\frac{1}{1-\alpha}}
\end{equation}
coming directly, where $\bar{\lambda}$ is the normalizing 
constant.
\medskip

Let us start with $M_{\alpha}(f)$ of (20) 
under the assumptions that $f(x)\ge 0$ for all $x$, 
$\int_a^bf(x){\rm d}x=1$, $\int_a^bx^{\delta}f(x){\rm d}x$ is 
fixed for all functional $f$ and for a specified $\delta>0$, 
$f(a)$ is the same for all functional $f$, $f(b)$ is the same 
for all functional $f$, for some limits $a$ and $b$, then the 
Euler equation becomes
\begin{equation}
(2-\alpha)f^{1-\alpha}-\lambda_1-\lambda_2x^{\delta}=0
\Rightarrow f=c_1[1+c_1^{*}x^{\delta}]^{\frac{1}{1-\alpha}}.
\end{equation}
If $c_1^{*}$ is written as $-s(1-\alpha), ~s>0$ then we have, 
writing $f_1$ for $f$,
\begin{equation}
f_1=c_1[1-s(1-\alpha)x^{\delta}]^{\frac{1}{1-\alpha}},~\delta>0,
\alpha<1,~0\le x\le \frac{1}{[s(1-\alpha)]^{\frac{1}{\delta}}}
\end{equation}
where $1-s(1-\alpha)x^{\delta}>0$. For $\alpha <1$ or 
$-\infty<\alpha<1$ the right side of (31) remains as a 
generalized type-1 beta model with the corresponding normalizing
constant $c_1$. For $\alpha >1$, writing $1-\alpha=-(\alpha-1)$
the model in (31) goes to a generalized type-2 beta form, 
namely,
\begin{equation}
f_2=c_2[1+s(\alpha-1)x^{\delta}]^{-\frac{1}{\alpha-1}}.
\end{equation}
When $\alpha\rightarrow 1$ in (31) or in (32) we have an 
extended or stretched exponential form,
\begin{equation}
f_3=c_3{\rm e}^{-s~x^{\delta}}.
\end{equation}
If $c_1^{*}$ in (30) is taken as positive then (30) for 
$\alpha<1,~\alpha>1,~\alpha\rightarrow 1$ will be increasing 
exponentially. Hence all possible forms are available from 
(30). The model in (31) is a special case of the distributional pathway 
model and for a discussion of the matrix-variate pathway model 
see Mathai (2005). Special cases of (31) and (32) for 
$\delta=1$ are Tsallis statistics (Gell-Mann and Tsallis, 2004; Ferri et al. 2005).\\
\medskip

Instead of optimizing $M_{\alpha}(f)$ of (22) under the 
conditions that $f(x)\ge 0$ for all $x$, $\int_a^bf(x){\rm d}x=1$ 
and $\int_a^bx^{\delta}f(x){\rm d}x$ is fixed, let us optimize 
under the following conditions: $f(x)\ge 0$ for all $x$, 
$\int_a^bf(x){\rm d}x<\infty$ and the following two moment-like 
expressions are fixed quantities for all functional $f$,

$$
\int_a^bx^{(\gamma-1)(1-\alpha)}f(x){\rm d}x=\mbox{ fixed  }, 
\int_a^bx^{(\gamma-1)(1-\alpha)+\delta}f(x){\rm d}x=\mbox{ fixed}.
$$
Then the Euler equation becomes
\begin{eqnarray*}
(2-\alpha)f^{1-\alpha}& &-\lambda_1x^{(\gamma-1)(1-\alpha)}
-\lambda_2x^{(\gamma-1)(1-\alpha)+\delta}=0\Rightarrow\\
f&=&c~x^{\gamma-1}[1+c^{*}x^{\delta}]^{\frac{1}{1-\alpha}}
\end{eqnarray*}
and for $c^{*}=-s(1-\alpha)$, $s>0$, we have the distributional pathway 
model for the real scalar case, namely
\begin{equation}
f(x)=c~x^{\gamma-1}[1-s(1-\alpha)x^{\delta}]^{\frac{1}{1-\alpha}},~\delta>0, ~s>0
\end{equation}
where $c$ is the normalizing constant. For $\alpha<1$, (34) 
gives a generalized type-1 beta form, for $\alpha>1$ it gives a 
generalized type-2 beta form and for $\alpha\rightarrow 1$ we 
have a generalized gamma form. For $\alpha>1$, (34) gives the 
superstatistics of Beck (2006) and Beck and Cohen (2003). For 
$\gamma=1, \delta=1$, (34) gives Tsallis statistics (Tsallis 2004, Cohen 2005).  
Densities appearing in a number of physical problems
are seen to be special cases of (34), a discussion of 
which may be seen from Mathai and Haubold (2006a). For example,
(34) for $\delta=2,\gamma=3,\alpha\rightarrow 1, x>0$ is the 
Maxwell-Boltzmann density; for $\delta=2,\gamma=1,\alpha
\rightarrow 1,-\infty<x<\infty$ is the Gaussian density; for 
$\gamma=\delta,\alpha\rightarrow 1$ is the Weibull density. 
For $\gamma=1,\delta=2,1<q<3$ we have the Wigner function $W(p)$ 
giving the atomic moment distribution in the framework of 
Fokker-Planck equation, see Douglas, Bergamini, and Renzoni (2006) where
\begin{equation}
W(p)=z_q^{-1}[1-\beta(1-q)p^2]^{\frac{1}{1-q}},~1<q<3.
\end{equation}
Before closing this section we may observe one more property 
for $M_{\alpha}(f)$. As an expected value
\begin{equation}
M_{\alpha}(f)=\frac{1}{\alpha-1}\left[E[f(x)]^{1-\alpha}-1\right].
\end{equation}
But Kerridge's (Kerridge, 1961) measure of ``inaccuracy" in 
assigning $q(x)$ for the true density $f(x)$, in the 
generalized form is 
\begin{equation}
H_{\alpha}(f:q)=\frac{1}{(2^{1-\alpha}-1)}\left[E[q(x)]^{\alpha-1}-1\right],
\end{equation}
which is also connected to the measure of directed divergence 
between $q(x)$ and $f(x)$. In (37) the normalizing constant 
is $2^{1-\alpha}-1$, the same factor appearing in Havrda-Charv\'t 
entropy. With different normalizing constants, as seen before,
(36) and (37) have the same forms as an expected value 
with $q(x)$ replaced by $f(x)$ in (36). Hence 
$M_{\alpha}(f)$ can also be looked upon as a type of directed 
divergence or ``inaccuracy" measure.
\medskip

\noindent
\textbf{3.\hskip.4cm Differential Equations}
\medskip

The functional part in (34), for a more general exponent, namely
\begin{equation}
g(x)=\frac{f(x)}{c}=x^{\gamma-1}[1-s(1-\alpha)x^{\delta}]^{\frac{\beta}{1-\alpha}},
~\alpha\ne 1, \delta>0, \beta>0, s>0
\end{equation}
is seen to satisfy the following differential equation for 
$\gamma\ne 1$ which defines the differential pathway.
\begin{eqnarray}
x\frac{{\rm d}}{{\rm d}x}g(x)&=&(\gamma-1)x^{\gamma-1}
[1-s(1-\alpha)x^{\delta}]^{\frac{\beta}{1-\alpha}}\nonumber\\
&& -s\beta\delta x^{\delta+\gamma-1}[1-s(1-\alpha)x^{\delta}]^{\frac{\beta}{1-\alpha}[1-\frac{(1-\alpha)}{\beta}]}.
\end{eqnarray}
Then for $\delta=\frac{(\gamma-1)(\alpha-1)}{\beta},
~\gamma\ne 1,~\alpha>1$ we have
\begin{eqnarray}
x\frac{{\rm d}}{{\rm d}x}g(x)&=&(\gamma-1)g(x)-s\beta\delta[g(x)]^{1-\frac{(1-\alpha)}{\beta}}\\
&=&(\gamma-1)g(x)-s\delta[g(x)]^{\alpha}\\
&&\mbox{  for  }\beta=1,\gamma\ne 1,\delta=(\gamma-1)(\alpha-1),\alpha>1.\nonumber
\end{eqnarray}
For $\gamma=1,\delta=1$ in (38) we have
\begin{eqnarray}
\frac{{\rm d}}{{\rm d}x}g(x)&=&-s[g(x)]^{\eta},~\eta=1-\frac{(1-\alpha)}{\beta}\\
&=&-s[g(x)]^{\alpha}\mbox{  for  }\beta=1.
\end{eqnarray}
Here (43) is the power law coming from Tsallis statistics (Gell-Mann and Tsallis, 2004).\\
\noindent
{\bf Acknowledgement} The authors would like to thank the Department of Science and Technology, Government of India, New Delhi, for the financial assistance for this work under project No. SR/S4/MS:287/05 which enabled this collaboration possible.\\
\medskip

\noindent
\textbf{4.\hskip.4cm References}
\medskip

\noindent
Beck, C. (2006). Stretched exponentials from superstatistics. 
{\it Physica A}, {\bf 365}, 96-101.\\

\noindent
Beck, C. and Cohen, E.G.D. (2003). Superstatistics. {\it Physica A},
{\bf 322}, 267-275.\\

\noindent
Cohen, E.G.D. (2005). Boltzmann and Einstein: Statistics and dynamics - An unsolved problem. 
{\it Pramana}, {\bf 64}, 635-643.\\

\noindent
Douglas, P., Bergamini, S., and Renzoni, F. (2006). Tunable 
Tsallis distribution in dissipative optical lattices. 
{\it Physical Review Letters}, {\bf 96}, 110601-1-4.\\

\noindent
Ferri, G.L., Martinez, S., and Plastino, A. (2005). Equivalence of the four versions of Tsallis's statistics.
{\it Journal of Statistical Mechanics: Theory and Experiment}, PO4009.\\

\noindent
Gell-Mann, M. and Tsallis, C. (Eds.) (2004). {\it Nonextensive Statistical Mechanics: 
Interdisciplinary Applications}. Oxford University Press, Oxford.\\

\noindent
Havrda, J. and Charv\'at, F. (1967). Quantification method of 
classification procedures: Concept of structural $\alpha$-entropy.
{\it Kybernetika}, {\bf 3}, 30-35.\\

\noindent
Kerridge, D.F. (1961). Inaccuracy and inference. 
{\it Journal of the Royal Statistical Society Series B}, 
{\bf 23}, 184-194.\\

\noindent
Mathai, A.M. (2005). A pathway to matrix-variate gamma and 
normal densities. {\it Linear Algebra and Its Applications}, 
{\bf 396}, 317-328.\\

\noindent
Mathai, A.M. and Haubold, H.J. (2006). Pathway model, Tsallis 
statistics, superstatistics and a generalized measure of 
entropy. {\it Physica A }, {\bf 375)}, 110-122.\\

\noindent
Mathai,A.M. and Haubold, H.J. (2006a). On generalized distributions and pathways. 
arXiv:cond-mat/0609526v2.\\

\noindent
Mathai, A.M. and Rathie, P.N. (1975). {\it Basic Concepts in 
Information Theory and Statistics: Axiomatic Foundations and 
Applications}, Wiley Halstead, New York and Wiley Eastern, New 
Delhi.\\

\noindent
R\'enyi, A. (1961). On measure of entropy and information. 
{\it Proceedings of the Fourth Berkeley Symposium on Mathematical 
Statistics and Probability, 1960}, University of California 
Press, 1961, Vol. 1, 547-561.\\

\noindent
Shannon, C.E. (1948). A mathematical theory of communication. 
{\it Bell System Technical Journal}, {\bf 27}, 379-423, 547-561.\\

\noindent
Tsallis, C. (1988). Possible generalization of Boltzmann-Gibbs 
statistics. {\it Journal of Statistical Physics}, {\bf 52}, 
479-487.\\

\noindent
Tsallis, C. (2004). What should a statistical mechanics satisfy
to reflect nature?, {\it Physica D}, {\bf 193}, 3-34.
\end{document}